\documentclass{amsart}
\usepackage{amssymb}
\usepackage{amsmath}
\usepackage{amsthm}
\newcommand{\pn}{\par\noindent}
\newcommand{\pmn}{\par\medskip\noindent}
\newtheorem{theor}{Theorem}[section]
\theoremstyle{definition} \newtheorem{defin}{Definition}[section]
\newtheorem{ex}{Example}[section] \theoremstyle{remark}
\newtheorem{rem}{Remark}[section]
\begin{document}
\title{On one property of Catalan numbers}
\author{Yury Kochetkov}
\date{}
\email{yukochetkov@hse.ru,yuyukochetkov@gmail.com}

\begin{abstract} We give a new proof of the following statement:
the Catalan number $C_n$ is divisible
by $n+2$, if $n$ is odd and $n\not\equiv 1\text{ mod }3$.
\end{abstract}

\maketitle

\section{Introduction}
\pn A plane tree is a tree imbedded in the plane. Two plane trees are the
same if there exists an orientation preserving homeomorphism that maps one
tree into another. An automorphism of a plane tree is a rotation around its
center of symmetry. By $\#\,{\rm Aut}(T)$ will be denoted the order of
the group of automorphisms of a tree $T$. \pmn
\begin{ex} For trees
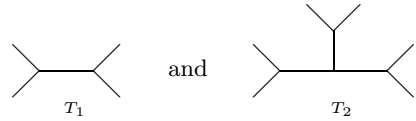
\[\begin{picture}(150,50)
\multiput(0,10)(30,10){2}{\line(1,1){10}}
\multiput(0,30)(30,-10){2}{\line(1,-1){10}}
\put(10,20){\line(1,0){20}} \put(19,4){\tiny $T_1$}

\put(58,18){\small and}

\multiput(90,10)(50,10){2}{\line(1,1){10}}
\multiput(90,30)(50,-10){2}{\line(1,-1){10}}
\put(100,20){\line(1,0){40}} \put(119,4){\tiny $T_2$}
\put(120,20){\line(0,1){15}} \put(120,35){\line(1,1){10}}
\put(120,35){\line(-1,1){10}}
\end{picture}\] we have: $\#\,{\rm Aut}(T_1)=2$ and $\#\,{\rm
Aut}(T_2)=3$.\end{ex}

\begin{defin} The number $\frac{1}{\#\,{\rm Aut}(T)}$ will be
called the \emph{weight} of a tree $T$ and will be denoted $w(T)$.
\end{defin}

\begin{defin} Let a tree $T$ has $a_1$ vertices of degree $1$,
$a_2$ vertices of degree $2$, $\ldots$ , and $a_k$ vertices of
degree $k$ ($k$ is the maximal degree). Then the array
$(a_1,a_2,\ldots,a_k)$ will be called the \emph{dual passport} of
$T$. \end{defin}

\begin{ex} At the figure below are presented all plane trees with
the dual passport $(4,2,0,1)$:
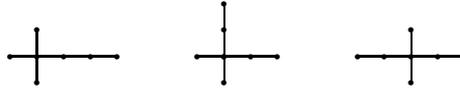
\[\begin{picture}(170,40) \put(0,15){\line(1,0){40}}
\put(10,5){\line(0,1){20}} \multiput(0,15)(10,0){5}{\circle*{2}}
\multiput(10,5)(0,20){2}{\circle*{2}}

\put(70,15){\line(1,0){30}} \put(80,5){\line(0,1){30}}
\multiput(70,15)(10,0){4}{\circle*{2}}
\multiput(80,5)(0,10){4}{\circle*{2}}

\put(130,15){\line(1,0){40}} \put(150,5){\line(0,1){20}}
\multiput(130,15)(10,0){5}{\circle*{2}}
\multiput(150,5)(0,20){2}{\circle*{2}}
\end{picture}\] Weights of these trees are $1$, $1$ and $\frac
12$, respectively. \end{ex} \pn We will use the following version
of the Goulden-Jackson theorem \cite{GJ}: \pmn {\bf Theorem.}
\emph{Let $S=\{T_1,\ldots,T_m\}$ be the set of all plane trees
with the dual passport $(a_1,\ldots,a_k)$, then}
$$\sum_{i=1}^m w(T_i)=\dfrac{(n-2)!}{a_1!\cdot
a_2!\cdot\ldots\cdot a_k!}\,,$$ \emph{where $n=a_1+\ldots+a_k$ is
the number of vertices.}

\begin{rem} For dual passport $(4,2,0,1)$ (Example 1.2) we have that $n=7$ and
$$1+1+\frac 12=\dfrac{5!}{4!\cdot 2!}\,.$$
\end{rem}

\section{Main result}
\pn Let us consider the set $S$ of all plane trees with the dual
passport $(n+2,0,n)$.

\begin{ex} For $n=4$ the set $S$ contains the following trees:
\[\begin{picture}(330,45) \put(0,15){\line(1,0){60}}
\multiput(10,15)(10,0){2}{\line(0,1){10}}
\multiput(40,15)(10,0){2}{\line(0,1){10}}

\put(80,15){\line(1,0){60}}
\multiput(90,15)(30,0){2}{\line(0,1){10}}
\multiput(100,15)(30,0){2}{\line(0,-1){10}}
\put(110,15){\circle*{3}}

\put(170,15){\line(1,0){60}}
\multiput(180,15)(30,0){2}{\line(0,-1){10}}
\multiput(190,15)(30,0){2}{\line(0,1){10}}
\put(200,15){\circle*{3}}

\multiput(260,5)(50,10){2}{\line(1,1){10}}
\multiput(260,25)(50,-10){2}{\line(1,-1){10}}
\put(270,15){\line(1,0){40}} \put(290,15){\line(0,1){15}}
\put(290,30){\line(1,1){10}} \put(290,30){\line(-1,1){10}}
\put(290,15){\circle*{3}}
\end{picture}\] with weights $1, \frac 12, \frac 12, \frac 13$,
respectively. Here black circles are centers of symmetry.\end{ex}

\begin{theor} The $n$-th Catalan number $C_n$ is divisible by
$n+2$, if $n$ is odd and $n\not\equiv 1\,{\rm mod}\,3$.\end{theor}

\begin{rem} For Catalan numbers see \cite{St}.\end{rem}

\begin{proof} If $T\in S$ is symmetric, then either its center of
symmetry is the middle of an edge, or a vertex of degree 3. In the
first case $T$ has an even number of vertices of degree 3. In the
second case $n=3p+1$ (one vertex of degree 3 is in the center,
others are symmetrically positioned around it). If there are no
symmetric trees in $S$, then all weights are equal to one. Thus,
the number
$$\frac{(2n)!}{(n+2)!\cdot n!}=\frac{1}{n+2}\cdot C_n$$ is an
integer. \end{proof}

\vspace{5mm}
\end{document}